\documentclass{amsproc}

\usepackage{amssymb}
\usepackage{graphicx}
\usepackage{amscd}
\usepackage{amsmath}
\theoremstyle{plain}

\newtheorem{definition}{Definition}

\newtheorem{lemma}{Lemma}

\newtheorem{remark}{Remark}

\newtheorem{theorem}{Theorem}

\numberwithin{equation}{section}

\begin{document}

\title{Dynamic inverse problem for complex Jacobi matrices.}

\author{ A. S. Mikhaylov}
\address{St. Petersburg   Department   of   V.A. Steklov    Institute   of   Mathematics
of   the   Russian   Academy   of   Sciences, 7, Fontanka, 191023
St. Petersburg, Russia and Saint Petersburg State University,
St.Petersburg State University, 7/9 Universitetskaya nab., St.
Petersburg, 199034 Russia.} \email{mikhaylov@pdmi.ras.ru}

\author{V. S. Mikhaylov}
\address{St.Petersburg   Department   of   V.A.Steklov    Institute   of   Mathematics
of   the   Russian   Academy   of   Sciences, 7, Fontanka, 191023
St. Petersburg, Russia.} \email{ftvsm78@gmail.com}


\maketitle






\noindent {\bf Abstract.} We consider the inverse dynamic problem
for a dynamical system with discrete time associated with a
semi-infinite complex Jacobi matrix. We propose two approaches of
recovering coefficients from dynamic response operator and answer
a question on the characterization of dynamic inverse data.

\section{Introduction.}

For a given sequence of complex numbers $\{a_0, a_1,\ldots\}$,
$\{b_1, b_2,\ldots  \}$, $a_i\not= 0$, we consider a dynamical
system with discrete time associated with a Jacobi marix:
\begin{equation}
\label{Jacobi_dyn} \left\{
\begin{array}l
u_{n,t+1}+u_{n,t-1}-a_{n}u_{n+1,t}-a_{n-1}u_{n-1,t}-b_nu_{n,t}=0,\quad n,t\in \mathbb{N},\\
u_{n,-1}=u_{n,0}=0,\quad n\in \mathbb{N}, \\
u_{0,t}=f_t,\quad t\in \mathbb{N}\cup\{0\},
\end{array}\right.
\end{equation}
which is a natural analog of a dynamical systems governed by a
wave equation on a semi-axis \cite{AM,BH,BM_1}. By an analogy with
continuous problems \cite{B07}, we treat the complex sequence
$f=(f_0,f_1,\ldots)$ as a \emph{boundary control}. The solution to
(\ref{Jacobi_dyn}) we denote by $u^f_{n,t}$. Having fixed $T\in
\mathbb{N}$, we associate the \emph{response operator} with
(\ref{Jacobi_dyn}), which maps the control $f=(f_0,\ldots
f_{T-1})$ to $u^f_{1,t}$:
\begin{equation*}
\left(R^T f\right)_t:=u^f_{1,t},\quad t=1,\ldots, T.
\end{equation*}
The inverse problem we will be dealing with consists in recovering
the sequences $\{b_1,b_2,\ldots,b_n\}$, $\{a_0,a_1,\ldots,a_n\}$
for some appropriated $n\in \mathbb{N}$ from $R^{T}$ with fixed
$T$. This problems is a natural discrete analog of an inverse
problem for a wave equation on a half-line where as an inverse
data the dynamic Dirichlet-to-Neumann map is used, see \cite{B07}.
Associated to this problem is a Jacobi matrix
\begin{equation} \label{Jac_matr}
A=\begin{pmatrix} b_1 & a_1 & 0 & 0 & 0 &\ldots \\
a_1 & b_2 & a_2 & 0 & 0 &\ldots \\
0 & a_2 & b_3 & a_3 & 0 & \ldots \\
\ldots &\ldots  &\ldots &\ldots & \ldots &\ldots
\end{pmatrix}.
\end{equation}
For $N\in \mathbb{N}$, by $A^N$ we denote the $N\times N$ Jacobi
matrix which is a block of (\ref{Jac_matr}) consisting of the
intersection of first $N$ columns with first $N$ rows of $A$.

We will use the nonselfadjoint variant \cite{AB,BH} of the
Boundary Control method \cite{B07} which was initially developed
for solving multidimensional dynamical inverse problems, but since
then was applied to multy- and one- dimensional inverse dynamical,
spectral and scattering problems, problems of signal processing
and identification problems. The application of the BC method to
one-dimensional problems are described in \cite{AM,BM_1}, the case
of real Jacobi matrix is considered in \cite{MM1,MM2,MM3}.



In the second section we study the forward problem: for
(\ref{Jacobi_dyn}) we prove the analog of Duhamel integral
representation formula. We also introduce the auxiliary problem,
introduce and derive representations for main operators of the
Boundary Control method: control and response operators for
(\ref{Jacobi_dyn}) and for auxiliary problem, and also derive the
representation for the connecting operator. In the third section
we outline two methods of recovering the unknown coefficients from
the response operator, namely Krein equations and factorization
method. Unlike the selfadjoint case, we will be able to recover
squares of $a_k,$ $k=1,2,\ldots$ only. We explain this feature of
the problem and answer a question on the characterization of
dynamic inverse data.

\section{Forward problem, auxiliary problem, operators of the Boundary Control method.}

We fix some positive integer $T$ and denote by $\mathcal{F}^T$ the
\emph{outer space} of the system (\ref{Jacobi_dyn}), the space of
controls: $\mathcal{F}^T:=\mathbb{C}^T$, $f\in \mathcal{F}^T$,
$f=(f_0,\ldots,f_{T-1})$, $f,g\in \mathcal{F}^T$,
$(f,g)_{\mathcal{F}^T}=\sum_{k=0}^{T-1} f_k\overline{g_k}$, we use
the notation $\mathcal{F}^\infty=\mathbb{C}^\infty$ when control
acts for all $t\geqslant 0$. We derive a representation formula
for the solution to (\ref{Jacobi_dyn}) which could be considered
as an analog of a Duhamel representation formula for an
initial-boundary value problem for a wave equation with a
potential on a half-line \cite{AM}.
\begin{lemma}
A solution to (\ref{Jacobi_dyn}) admits the representation
\begin{equation}
\label{Jac_sol_rep} u^f_{n,t}=\prod_{k=0}^{n-1}
a_kf_{t-n}+\sum_{s=n}^{t-1}w_{n,s}f_{t-s-1},\quad n,t\in
\mathbb{N},
\end{equation}
where $w_{n,s}$ satisfies the Goursat problem
\begin{equation}
\label{Goursat} \left\{
\begin{array}l
w_{n,s+1}+w_{n,s-1}-a_nw_{n+1,s}-a_{n-1}w_{n-1,s}-b_nw_{n,s}=\\
=-\delta_{s,n}(1-a_n^2)\prod_{k=0}^{n-1}a_k,\,n,s\in \mathbb{N}, \,\,s>n,\\
w_{n,n}-b_n\prod_{k=0}^{n-1}a_k-a_{n-1}w_{n-1,n-1}=0,\quad n\in \mathbb{N},\\
w_{0,t}=0,\quad t\in \mathbb{N}_0.
\end{array}
\right.
\end{equation}
\end{lemma}
\begin{proof}
We assume that $u^f_{n,t}$ has a form (\ref{Jac_sol_rep}) with
unknown $w_{n,s}$ and plug it into (\ref{Jacobi_dyn}):
\begin{eqnarray*}
0=\prod_{k=0}^{n-1}a_kf_{t+1-n}+\prod_{k=0}^{n-1}a_kf_{t-1-n}-a_n\prod_{k=0}^{n}a_kf_{t-n-1}-a_{n-1}\prod_{k=0}^{n-2}a_kf_{t-n+1}\\
-b_n\prod_{k=0}^{n-1}a_kf_{t-n}+\sum_{s=n}^{t}w_{n,s}f_{t-s}+\sum_{s=n}^{t-2}w_{n,s}f_{t-s-2}\\-a_n\sum_{s=n+1}^{t-1}w_{n+1,s}f_{t-s-1}
-a_{n-1}\sum_{s=n-1}^{t-1}w_{n-1,s}f_{t-s-1}-\sum_{s=n}^{t-1}b_nw_{n,s}f_{t-s-1}.
\end{eqnarray*}
Changing the order of summation and evaluating we get
\begin{eqnarray*}
0=\left(1-a_n^2\right)\prod_{k=0}^{n-1}a_kf_{t-n-1}-b_n\prod_{k=0}^{n-1}a_kf_{t-n}\\
-\sum_{s=n}^{t-1}f_{t-s-1}\left(b_nw_{n,s}+a_nw_{n+1,s}+a_{n-1}w_{n-1,s}\right)
+a_nw_{n+1,n}f_{t-n-1}\\
-a_{n-1}w_{n-1,n-1}f_{t-n}
+\sum_{s=n-1}^{t-1}w_{n,s+1}f_{t-s-1}+\sum_{s=n+1}^{t-1}w_{n,s-1}f_{t-s-1}\\
=\sum_{s=n}^{t-1}f_{t-s-1}\left(w_{n,s+1}+w_{n,s-1}-a_nw_{n+1,s}-a_{n-1}w_{n-1,s}-b_nw_{n,s}\right)\\
-b_n\prod_{k=0}^{n-1}a_kf_{t-n}+\left(1-a_n^2\right)\prod_{k=0}^{n-1}a_kf_{t-n-1}+a_nw_{n+1,n}f_{t-n-1}\\
-a_{n-1}w_{n-1,n-1}f_{t-n}+w_{n,n}f_{t-n}-w_{n,n-1}f_{t-n-1}.
\end{eqnarray*}
Finally we arrive at
\begin{eqnarray*}
\sum_{s=n}^{t-1}f_{t-s-1}\left(w_{n,s+1}+w_{n,s-1}-a_nw_{n+1,s}-a_{n-1}w_{n-1,s}-b_nw_{n,s}+\left(1-a_n^2\right)\prod_{k=0}^{n-1}a_k\delta_{sn}\right)\\
+f_{t-n}\left(w_{n,n}-a_{n-1}w_{n-1,n-1}-b_n\prod_{k=0}^{n-1}a_k\right)=0.
\end{eqnarray*}
Counting that $w_{n,s}=0$ when $n>s$ and arbitrariness of $f$, we
arrive at (\ref{Goursat}).
\end{proof}

\begin{definition}
For $f,g\in \mathcal{F}^\infty$ we define the convolution
$c=f*g\in \mathcal{F}^\infty$ by the formula
\begin{equation*}
c_t=\sum_{s=0}^{t}f_sg_{t-s},\quad t\in \mathbb{N}\cup \{0\}.
\end{equation*}
\end{definition}

We introduce the analog of the dynamic \emph{response operator}
(dynamic Dirichlet-to-Neumann map) \cite{B07} for the system
(\ref{Jacobi_dyn}):
\begin{definition}
For (\ref{Jacobi_dyn}) the \emph{response operator}
$R^T:\mathcal{F}^T\mapsto \mathbb{C}^T$ is defined by the rule
\begin{equation*}
\left(R^Tf\right)_t=u^f_{1,t}, \quad t=1,\ldots,T.
\end{equation*}
\end{definition}
The \emph{response vector} is a convolution kernel of a response
operator, $r=(r_0,r_1,\ldots,r_{T-1})=(a_0,w_{1,1},w_{1,2},\ldots
w_{1,T-1})$, in accordance with (\ref{Jac_sol_rep}):
\begin{eqnarray}
\label{R_def}
\left(R^Tf\right)_t=u^f_{1,t}=a_0f_{t-1}+\sum_{s=1}^{t-1}
w_{1,s}f_{t-1-s}
\quad t=1,\ldots,T.\\
\notag \left(R^Tf\right)=r*f_{\cdot-1}.
\end{eqnarray}
 By choosing the special control
 $f=\delta=(1,0,0,\ldots)$, the kernel of a response operator can be determined as
\begin{equation}
\label{con11} \left(R^T\delta\right)_t=u^\delta_{1,t}=
r_{t-1},\quad t=1,2,\ldots.
\end{equation}
The inverse problem we will be dealing with consists of recovering
a Jacobi matrix (i.e. the sequences $\{a_1,a_2,\ldots\}$, $\{b_1,
b_2,\ldots \}$) and $a_0$ from the response operator.

We introduce the \emph{inner space} of dynamical system
(\ref{Jacobi_dyn}) $\mathcal{H}^T:=\mathbb{C}^T$, $h\in
\mathcal{H}^T$, $h=(h_1,\ldots, h_T)$ with the inner product
$h,l\in \mathcal{H}^T$, $(h,g)_{\mathcal{H}^T}=\sum_{k=1}^T
h_k\overline{g_k}$. The \emph{control operator}
$W^T:\mathcal{F}^T\mapsto \mathcal{H}^T$ is defined by the rule
\begin{equation*}
W^Tf:=u^f_{n,T},\quad n=1,\ldots,T.
\end{equation*}
From (\ref{Jac_sol_rep}) we deduce the representation for $W^T$:
\begin{equation}
\label{W^T_rep} \left(W^Tf\right)_n=u^f_{n,T}=\prod_{k=0}^{n-1}
a_kf_{T-n}+\sum_{s=n}^{T-1}w_{n,s}f_{T-s-1},\quad n=1,\ldots,T.
\end{equation}
The following statement is equivalent to a boundary
controllability of (\ref{Jacobi_dyn}).
\begin{lemma}
\label{teor_control} The operator $W^T$ is an isomorphism between
$\mathcal{F}^T$ and $\mathcal{H}^T$.
\end{lemma}
\begin{proof}
We fix some $a\in \mathcal{H}^T$ and look for a control $f\in
\mathcal{F}^T$ such that $W^Tf=a$. We write down the action of the
operator $W^T$ as
\begin{equation}
\label{WT}
W^Tf=\begin{pmatrix} u_{1,T}\\
u_{2,T}\\
\cdot\\
u_{k,T}\\
\cdot\\
u_{T,T}
\end{pmatrix}=\begin{pmatrix}
a_0& w_{1,1} & w_{1,2} & \ldots & \ldots & w_{1,T-1}\\
0 & a_0a_1 & w_{2,2} & \ldots & \ldots & w_{2, T-1}\\
\cdot & \cdot & \cdot & \cdot & \cdot & \cdot \\
0 & \ldots & \prod_{j=0}^{k-1}
a_j& w_{k,k} & \ldots & w_{k, T-1}\\
\cdot & \cdot & \cdot & \cdot & \cdot & \cdot \\
0 & 0 & 0 & 0 & \ldots & \prod_{k=0}^{T-1} a_{T-1}
\end{pmatrix}
\begin{pmatrix} f_{T-1}\\
f_{T-2}\\
\cdot\\
f_{T-k-1}\\
\cdot\\
f_{0}
\end{pmatrix}.
\end{equation}
On introducing  the notations
\begin{eqnarray*}
J_T: \mathcal{F}^T\mapsto \mathcal{F}^T,\quad
\left(J_Tf\right)_n=f_{T-1-n},\quad n=0,\ldots, T-1,\\
A\in \mathbb{R}^{T\times T},\quad
a_{ii}=\prod_{k=0}^{i-1}a_i,\quad a_{ij}=0,\, i\not=j,\\
\\ K\in \mathbb{R}^{T\times T},\quad k_{ij}=0,\, i\geqslant j,
\, k_{ij}=w_{ij-1},\,i<j,
\end{eqnarray*}
we have that
\begin{equation}
\label{WT1} W^T=V^TJ^T=\left(A+K\right)J^T.
\end{equation}
Obviously, this operator is invertible, which completes the
statement of the lemma.
\end{proof}

Along with the system (\ref{Jacobi_dyn}) we consider the auxiliary
system associated with the complex conjugate matrix $\overline A$:
\begin{equation}
\label{Jacobi_dyn_aux}
\begin{cases}
v_{n,t+1}+v_{n,t-1}-\overline{a_{n}}v_{n+1,t}-\overline{a_{n-1}}v_{n-1,t}-\overline{b_n}v_{n,t}=0,\quad n,t\in \mathbb{N},\\
v_{n,-1}=v_{n,0}=0,\quad n\in \mathbb{N}, \\
v_{0,t}=f_t,\quad t\in \mathbb{N}\cup\{0\}.
\end{cases}
\end{equation}
The objects corresponding to the system (\ref{Jacobi_dyn_aux}) we
equip with the symbol $\#$. The direct calculations shows:
\begin{lemma}
The control and response operators of the system $\#$ a related
with control and response operators of the original system by
\begin{equation}
\label{adj_prop}
W^T_{\#}=\overline{W^T},\quad
R^T_{\#}=\overline{R^T},
\end{equation}
that is the matrix of $W^T_{\#}$ and response vector $r_{\#}$ are
complex conjugate to the matrix of $W^T$ and vector $r$.
\end{lemma}

For the systems (\ref{Jacobi_dyn}), (\ref{Jacobi_dyn_aux}) we
introduce the \emph{connecting operator} $C^T:
\mathcal{F}^T\mapsto \mathcal{F}^T$ by the quadratic form: for
arbitrary $f,g\in \mathcal{F}^T$ we define
\begin{equation}
\label{C_T_def} \left(C^T
f,g\right)_{\mathcal{F}^T}=\left(u^f_{\cdot,T},
v^g_{\cdot,T}\right)_{\mathcal{H}^T}=\left(W^Tf,W^T_\#g\right)_{\mathcal{H}^T}.
\end{equation}
The next statement will be very important in solving the dynamic
inverse problem:
\begin{theorem}
The connecting operator $C^T$ is an isomorphism in
$\mathcal{F}^T$, it admits the representation in terms of inverse
data:
\begin{equation}
\label{C_T_repr} C^T=a_0C^T_{ij},\quad
C^T_{ij}=\sum_{k=0}^{T-\max{i,j}}r_{|i-j|+2k},\quad r_0=a_0,
\end{equation}
\begin{equation*}
C^T=
\begin{pmatrix}
r_0+r_2+\ldots+r_{2T-2} & r_1+r_3+\ldots+r_{2T-3} & \cdot &
r_T+r_{T-2} &
r_{T-1}\\
r_1+r_3+\ldots+r_{2T-3} & r_0+r_2+\ldots+r_{2T-4} & \cdot & \ldots
&r_{T-2}\\
\cdot & \cdot & \cdot & \cdot & \cdot \\
r_{T-3}+r_{T-1}+r_{T+1} &\ldots & \cdot & r_1+r_3 & r_2\\
r_{T}+r_{T-2}&\ldots & \cdot &r_0+r_2&r_1 \\
r_{T-1}& r_{T-2}& \cdot & r_1 &r_0
\end{pmatrix}.
\end{equation*}
\end{theorem}
\begin{proof}
We observe that $C^T=\left(W^T_\#\right)^*W^T$, so due to Lemma
\ref{teor_control}, $C^T$ is an isomorphism in $\mathcal{F}^T$. We
use the variant of the BC method for  nonselfadjoint problems
\cite{AB}, for fixed $f,g\in \mathcal{F}^T$ we introduce the
\emph{Blagoveshchenskii function} by the rule
\begin{equation*}
\psi_{n,t}:=\left(u^f_{\cdot,n},
v^g_{\cdot,t}\right)_{\mathcal{H}^T}=\sum_{k=1}^T
u^f_{k,n}\overline{v^g_{k,t}}.
\end{equation*}
We show that $\psi_{n,t}$ satisfies certain difference equation.
Indeed, we can evaluate:
\begin{eqnarray*}
\psi_{n,t+1}+\psi_{n,t-1}-\psi_{n+1,t}-\psi_{n-1,t}=\sum_{k=1}^T
u^f_{k,n}\overline{\left(v^g_{k,t+1}+v^g_{k,t-1}\right)}\\
-\sum_{k=1}^T
\left(u^f_{k,n+1}+u^f_{k,n-1}\right)\overline{v^g_{k,t}}=
\sum_{k=1}^T
u^f_{k,n}\overline{\left(\overline{a_k}v^g_{k+1,t}+\overline{a_{k-1}}v^g_{k-1,t}+\overline{b_k}v^g_{k,t}\right)}
\\
-\sum_{k=1}^T
\left(a_ku^f_{k+1,n}+a_{k-1}u^f_{k-1,n}+b_ku^f_{k,n}\right)\overline{v^g_{k,t}}=\\
\sum_{k=1}^T
u^f_{k,n}{\left({a_k}\overline{v^g_{k+1,t}}+{a_{k-1}}\overline{v^g_{k-1,t}}\right)}-
\sum_{k=1}^T\left(a_ku^f_{k+1,n}+a_{k-1}u^f_{k-1,n}\right)\overline{v^g_{k,t}}=\\
+a_0u^f_{1,n}\overline{v^g_{0,t}}-a_0u^f_{0,n}\overline{v^g_{1,t}}
+a_Tu^f_{T,n}\overline{v^g_{T+1,t}}-a_Tu^f_{T+1,n}\overline{v^g_{T,t}}\\
=a_0\left[(Rf)_n\overline{g_t}-f_n\overline{(R_\#g)_t}\right].
\end{eqnarray*}
Thus we arrived at the following difference equation on
$\psi_{n,t}$:
\begin{eqnarray}
\label{Blag_eqn}
&\left\{
\begin{array}l
\psi_{n,t+1}+\psi_{n,t-1}-\psi_{n+1,t}-\psi_{n-1,t}=h_{n,t},\quad n,t\in \mathbb{N}_0,\\
\psi_{0,t}=0,\,\, \psi_{n,0}=0,
\end{array}
\right. \\
&h_{n,t}=a_0\left[g_t(Rf)_n-f_n(Rg)_t\right].\notag
\end{eqnarray}
We introduce the set
\begin{eqnarray*}
K(n,t):=\left\{(n,t)\cup \{(n-1,t-1),
(n+1,t-1)\}\cup\{(n-2,t-2),(n,t-2),\right.\\
\left.(n+2,t-2)\}\cup\ldots\cup\{(n-t,0),(n-t+2,0),\ldots,(n+t-2,0),(n+t,0)\}\right\}\\
=\bigcup_{\tau=0}^t\bigcup_{k=0}^\tau
\left(n-\tau+2k,t-\tau\right).
\end{eqnarray*}
The solution to (\ref{Blag_eqn}) is given by (see \cite{MM1,MM2})
\begin{equation*}
\psi_{n,t}=\sum_{(k,\tau)\in K(n,t-1)}h(k,\tau).
\end{equation*}
We observe that $\psi_{T,T}=\left(C^Tf,g\right)$, so
\begin{equation}
\label{C_T_sol}\left(C^Tf,g\right)=\sum_{(k,\tau)\in
K(T,T-1)}h(k,\tau).
\end{equation}
We note that in the right hand side of (\ref{C_T_sol}) the
argument $k$ runs from $1$ to $2T-1.$ Extending $f\in
\mathcal{F}^T$, $f=(f_0,\ldots,f_{T-1})$ to $f\in
\mathcal{F}^{2T}$ by:
\begin{equation*}
f_T=0,\quad f_{T+k}=-f_{T-k},\,\, k=1,2,\ldots, T-1,
\end{equation*}
we deduce that $\sum_{k,\tau\in K(T,T-1)} f_k(R^Tg)_\tau=0$, so
(\ref{C_T_sol}) leads to
\begin{eqnarray*}
\left(C^Tf,g\right)=\sum_{k,\tau\in K(T,T-1)}
g_\tau\left(R^{2T}f\right)_k=g_0\left[\left(R^{2T}f\right)_1+\left(R^{2T}f\right)_3+\ldots+\left(R^{2T}f\right)_{2T-1}\right]\\
+g_1\left[\left(R^{2T}f\right)_2+\left(R^{2T}f\right)_4+\ldots+\left(R^{2T}f\right)_{2T-2}\right]+\ldots+g_{T-1}\left(R^{2T}f\right)_{T}.
\end{eqnarray*}
The latter relation yields
\begin{equation*}
C^Tf=\left(\left(R^{2T}f\right)_1+\ldots+\left(R^{2T}f\right)_{2T-1},\left(R^{2T}f\right)_2+\ldots+\left(R^{2T}f\right)_{2T-2},\ldots,\left(R^{2T}f\right)_{T}
\right)
\end{equation*}
from where the statement of the theorem follows.
\end{proof}

By $B^\tau$ we denote the matrix transposed to $B\in
\mathbb{C}^{n\times n}$.

The relations (\ref{adj_prop}) imply the following
\begin{remark}
The connecting operator is complex symmetric:
\begin{equation*}
\left(C^T\right)^*=\overline{C^T}, \quad \text{or} \quad
\left(C^T\right)^\tau=C^T.
\end{equation*}
\end{remark}

\section{Inverse problem. }

Due to the final speed of wave propagation in ((\ref{Jacobi_dyn}))
the solution $u^f$ depends on the coefficients $a_n,b_n$ in the
following way: for $M\in \mathbb{N}$, $u^f_{M-1,M}$ depends on
$\{a_0,\ldots,a_{M-1}\}$, $\{b_1,\ldots,b_{M-1}\}$, which implies
that $u^f_{1,2M-1}$ depends of the same set of parameters:
\begin{remark}
\label{Rem1} The response $R^{2T}$ (or, what is equivalent, the
response vector $(r_0,r_1,\ldots,r_{2T-2})$) depends on
$\{a_0,\ldots,a_{T-1}\}$, $\{b_1,\ldots,b_T\}$.
\end{remark}
Thus the natural set up of the dynamic inverse problem for
(\ref{Jacobi_dyn}): by the given operator $R^{2T}$ to recover
$\{a_0,\ldots,a_{T-1}\}$ and $\{b_1,\ldots,b_{T-1}\}$.

We also note that $a_0=r_0$, which follows from \ref{R_def}).

\subsection{Krein equations}
Let $\alpha,\beta\in \mathbb{R}$ and $y$ be a solution to
\begin{equation}
\label{y_special} \left\{
\begin{array}l
a_ky_{k+1}+a_{k-1}y_{k-1}+b_ky_k=0,\\
y_0=\alpha,\,\, y_1=\beta.
\end{array}
\right.
\end{equation}
We set up the following control problem: to find a control $f^T\in
\mathcal{F}^T$ such that
\begin{equation}
\label{Control_probl}
\left(W^Tf^T\right)_k=y_k,\quad
k=1,\ldots,T.
\end{equation}
Due to Lemma \ref{teor_control}, this control problem has a unique
solution. Let $\varkappa^T$ be a solution to
\begin{equation}
\label{kappa} \left\{
\begin{array}l
\varkappa^T_{t+1}+\varkappa^T_{t-1}=0,\quad t=0,\ldots,T,\\
\varkappa^T_{T}=0,\,\, \varkappa^T_{T-1}=1.
\end{array}
\right.
\end{equation}
We show that the control $f^T$ satisfies the Krein equation:
\begin{theorem}
The control $f^T$, defined by (\ref{Control_probl}) satisfies the
following equation in $\mathcal{F}^T$:
\begin{equation}
\label{C_T_Krein} C^Tf^T=a_0\left[\beta\varkappa^T-\alpha
\left(R^T_{\#}\right)^*\varkappa^T\right].
\end{equation}
\end{theorem}
\begin{proof}
Let $f^T$ be a solution to (\ref{Control_probl}). We observe that
for any fixed $g\in \mathcal{F}^T$:
\begin{equation}
\label{Kr_1}
v^g_{k,T}=\sum_{t=1}^{T-1}\left(v^g_{k,t+1}+v^g_{k,t-1}\right)\varkappa^T_t,\quad
k\leqslant T.
\end{equation}
Indeed, changing the order of summation in the right hand side of
(\ref{Kr_1}) yields
\begin{equation*}
\sum_{t=1}^{T-1}\left(v^g_{k,t+1}+v^g_{k,t-1}\right)\varkappa^T_t=\sum_{t=1}^{T-1}\left(\varkappa^T_{t+1}+\varkappa^T_{t-1}
\right)v^g_{k,t}+v^g_{k,0}\varkappa^T_1-v^g_{k,T}\varkappa^T_{T-1},
\end{equation*}
which gives (\ref{Kr_1}) due to (\ref{kappa}). Using this
observation, we can evaluate
\begin{eqnarray*}
\left(C^Tf^T,g\right)=\sum_{k=1}^T
y_k\overline{v^g_{k,T}}=\sum_{k=1}^T\sum_{t=0}^{T-1}\varkappa^T_t
y_k\overline{\left(v^g_{k,t+1}+v^g_{k,t-1}\right)}\\
=\sum_{t=0}^{T-1}\varkappa^T_t\left(\sum_{k=1}^T
\left(a_k\overline{v^g_{k+1,t}}y_k+a_{k-1}\overline{v^g_{k-1,t}}y_k + b_k\overline{v^g_{k,t}}y_k\right)\right)\\
=\sum_{t=0}^{T-1}\varkappa^T_t\left(\sum_{k=1}^T
\left(\overline{v^g_{k,t}}(a_ky_{k+1}+a_{k-1}y_{k-1}+b_ky_k\right)+a_0\overline{v^g_{0,t}}y_1\right.\\
\left.+a_T\overline{v^g_{T+1,t}}y_T-a_0\overline{v^g_{1,t}}y_0-a_T\overline{v^g_{T,t}}y_{T+1}
\right) =\sum_{t=0}^{T-1}\varkappa^T_t\left(a_0\beta
\overline{g_t}-a_0\alpha\overline{\left(R^T_{\#}g\right)_t} \right)\\
=\left(\varkappa^T, \overline{a_0}\left[\overline{\beta}
g-\overline{\alpha}
\left(R^T_{\#}g\right)\right]\right)_{\mathcal{F}^T}=\left(a_0\left[\beta\varkappa^T
- {\alpha} \left(\left(R^T_{\#}\right)^*\varkappa^T\right)\right],
g\right)_{\mathcal{F}^T}.
\end{eqnarray*}
From where (\ref{C_T_Krein}) follows.
\end{proof}

Now we describe the procedure of the recovering $a_0,$ $a_n,$
$b_n$, $n=1,\ldots,T-1$ from the solutions of Krein
(\ref{C_T_Krein}) equations $f^\tau\in \mathcal{F}^\tau$ for
$\tau=1,\ldots,T$. From (\ref{Jac_sol_rep}) and (\ref{Goursat}) we
infer that
\begin{eqnarray*}
u^{f^T}_{T,T}=\prod_{k=0}^{T-1}a_k f^T_0,\\
u^{f^T}_{T-1,T}=\prod_{k=0}^{T-2}a_k
f^T_1+\prod_{k=0}^{T-2}a_k(b_1+b_2+\ldots+b_{T-1})f^T_0.
\end{eqnarray*}
Notice that we know $a_0=r_0$. Let $T=2$, then we have:
\begin{eqnarray}
y_2=u^{f^2}_{2,2}=a_0a_1f^2_0,\label{y_2}\\
y_1=u^{f^2}_{1,2}=a_0f^2_1+a_0b_1f^2_0, \label{y_1}
\end{eqnarray}
In (\ref{y_1}) we know $y_1=\beta,$ $a_0,$ $f^2_1,$ $f^2_0$, so we
can recover $b_1$. On the other hand, using (\ref{y_special}), we
have a system
\begin{equation*}
\left\{
\begin{array}l
y_2=a_0a_1f^2_0,\\
a_1y_2+a_0\alpha+b_1\beta=0
\end{array}
\right.
\end{equation*}
Since $a_1\not=0,$ we can recover $(a_1)^2.$ We proceed by the
induction: assuming that we have already recovered $b_{k-2},$
$\left((a_{k-2}\right)^2$ for $k\leqslant n$ and we know that
$y_{k-1}=\prod_{i=0}^{k-2}a_if^{k-1}_0$, we recover
$\left(a_{n-1}\right)^2,$ $b_{n-1}$. Bearing in mind that
\begin{eqnarray}
y_n=u^{f^n}_{n,n}=\prod_{k=0}^{n-2}a_k a_{n-1}f^n_0,\label{y_n}\\
y_{n-1}=u^{f^n}_{n-1,n}=\prod_{k=0}^{n-2}a_kf^{n}_1+\prod_{k=0}^{n-2}a_k(b_1+\ldots+b_{n-2}+b_{n-1})f^n_0,
\label{y_n1}
\end{eqnarray}
and that we know $f^n_0$, $f^n_1$, and $(a_k)^2,$ $b_k$,
$k\leqslant n-2,$ we plug $y_{n-1}=\prod_{k=0}^{n-2}a_kf^{n-1}_0$
into (\ref{y_n1}), cancel out $\prod_{k=0}^{n-2}a_k$ and recover
$b_{n-1}$. Using (\ref{y_special}) and (\ref{y_n}) leads to the
following relations
\begin{eqnarray*}
y_n=\prod_{k=0}^{n-2}a_k a_{n-1}f^n_0,\\
a_{n-1}y_n+a_{n-2}y_{n-2}+b_{n-1}y_{n-1}=0.
\end{eqnarray*}
We rewrite the latter one using representations for $y_n$,
$y_{n-1}$ and $y_{n-2}$:
\begin{equation*}
\left(a_{n-1}\right)^2\prod_{k=0}^{n-2}a_k f_0^n+
\prod_{k=0}^{n-2}a_k f_0^{n-2}+b_{n-1}\prod_{k=0}^{n-2}a_k
f_0^{n-1}=0,
\end{equation*}
from where we can recover $\left(a_{n-1}\right)^2$.

\begin{remark}
We see that the described procedure allows one to recover
$\left(a_k\right)^2$, $k=1,\ldots$ only.
\end{remark}

\subsection{Factorization method}

We make use of the fact that the matrix $C^T$ has a special
structure (\ref{C_T_def}) -- it is a product of a triangular
matrix and its transposed. We rewrite the operator $W^T$ as $W^T=
V^TJ_T$ where
\begin{equation*}
W^Tf=\begin{pmatrix}
a_0 & w_{1,1} & w_{1,2} & \ldots & w_{1,T-1}\\
0 & a_0a_1 & w_{2,2} &  \ldots & w_{2, T-1}\\
\cdot & \cdot & \cdot & \cdot & \cdot \\
0 & \ldots & \prod_{j=1}^{k-1}a_j& \ldots & w_{k, T-1}\\
\cdot & \cdot & \cdot  & \cdot & \cdot \\
0 & 0 & 0  & \ldots & \prod_{j=1}^{T-1}a_j
\end{pmatrix}
\begin{pmatrix}
0 & 0 & 0 & \ldots & 1\\
0 & 0 & 0 & \ldots  & 0\\
\cdot & \cdot & \cdot & \cdot &  \cdot \\
0 & \ldots & 1& 0  & 0\\
\cdot & \cdot & \cdot & \cdot &  \cdot \\
1 & 0 & 0 & 0 &  0
\end{pmatrix}
\begin{pmatrix} f_{0}\\
f_{2}\\
\cdot\\
f_{T-k-1}\\
\cdot\\
f_{T-1}
\end{pmatrix}
\end{equation*}
Using the definition (\ref{C_T_def}) and the invertibility of
$W^T$ (cf. Lemma \ref{teor_control}) gives that
\begin{equation*}
C^T=\left(W^T_{\#}\right)^*W^T,\quad \text{or} \quad
\left(\left(W^T_{\#}\right)^{-1}\right)^*C^T\left(W^T\right)^{-1}=I.
\end{equation*}
Using (\ref{adj_prop}) we can rewrite the latter equality as
\begin{equation}
\label{C_T_eqn_ker} \left(\left( V^T\right)^{-1}\right)^\tau
C_T\left( V^T\right)^{-1}=I,\quad C_T=J_TC^TJ_T,
\end{equation}
where the matrix $C_T$ has the entries:
\begin{equation}
\label{C_overline_repr}
c_{i,j}=\left\{C_T\right\}_{ij}=C_{T+1-j,T+1-i},\quad C_T=a_0
\begin{pmatrix}
r_0 & r_1 & r_2 & \ldots & r_{T-1}\\
r_1 & r_0+r_2 & r_1+r_3 & \ldots
&..\\
r_2 & r_1+r_3 & r_0+r_2+r_4 & \ldots & ..\\
\cdot & \cdot & \cdot & \cdot & \cdot \\
\end{pmatrix},
\end{equation}
and  $\left(V^T\right)^{-1}$ has a form
\begin{equation}
\label{W_T_bar} \left(V^T\right)^{-1}=\begin{pmatrix}
a_{1,1} & a_{1,2} & a_{1,3}& \ldots & a_{1,T} \\
0 & a_{2,2} & a_{2,3} &\ldots &..\\
\cdot & \cdot & \cdot & a_{T-1,T-1} & a_{T-1,T} \\
0 &\ldots &\ldots & 0 & a_{T,T}
\end{pmatrix}.
\end{equation}
We multiply the $k-$th row of $V^T$ by $k-$th column of
$\left(V^T\right)^{-1}$ to get $a_{k,k}a_0a_1\ldots a_{k-1}=1$, so
diagonal elements of (\ref{W_T_bar}) satisfy the relation:
\begin{equation}
\label{AK_form} a_{k,k}=\left(\prod_{j=0}^{k-1}a_j\right)^{-1}.
\end{equation}
Multiplying the $k-$th row of $V^T$ by $k+1-$th column of
$\left(V^T\right)^{-1}$ leads to the relation
\begin{equation*}
a_{k,k+1}a_0a_1\ldots a_{k-1}+a_{k+1,k+1}w_{k,k}=0,
\end{equation*}
from where we deduce that
\begin{equation}
\label{AK1_form}
a_{k,k+1}=-\left(\prod_{j=0}^{k}a_j\right)^{-2}a_k w_{k,k},\quad
k=1,\ldots, T-1.
\end{equation}

All aforesaid leads to the equivalent form of (\ref{C_T_eqn_ker}):
\begin{equation}
\label{matr_eq}
\begin{pmatrix}
{a_{1,1}} & 0 & . &  0 \\
{a_{1,2}} & {a_{2,2}} & 0  &.\\
\cdot & \cdot & \cdot & \cdot  \\
{a_{1,T}} & . & .  & {a_{T,T}}
\end{pmatrix}
\begin{pmatrix}
 c_{11} & .. & .. &   c_{1T} \\
.. & .. & ..  &..\\
\cdot & \cdot & \cdot & \cdot  \\
 c_{T1} &.. &   &  c_{TT}
\end{pmatrix}\begin{pmatrix}
a_{1,1} & a_{1,2} & ..& a_{1,T} \\
0 & a_{2,2} & ..  & a_{2,T}\\
\cdot & \cdot & \cdot & \cdot  \\
0 &\ldots &\ldots  & a_{T,T}
\end{pmatrix}=I
\end{equation}
In the above equality $c_{ij}$ are given (see
(\ref{C_overline_repr})), the entries $a_{ij}$ are unknown. A
direct consequence of (\ref{matr_eq}) is an equality for
determinants:
\begin{equation*}
\det{\left(\left( V^T\right)^{-1}\right)^\tau}\det{
C_T}\det{\left(\left( V^T\right)^{-1}\right)}=1,
\end{equation*}
which yields
\begin{equation*}
\left(a_{1,1}\right)^2*\ldots*(a_{T,T})^2=\left(\det{
C_T}\right)^{-1}.
\end{equation*}
From the above equality we derive that
\begin{equation*}
(a_{1,1})^2=\left(\det{C_1}\right)^{-1},\quad
(a_{2,2})^2=\left(\frac{\det{C_2}}{\det{ C_1}}\right)^{-1},\quad
(a_{k,k})^2=\left(\frac{\det{C_k}}{\det{ C_{k-1}}}\right)^{-1}.
\end{equation*}
Combining latter relations with (\ref{AK_form}), we deduce that
\begin{equation*}
(a_0)^2*\ldots*(a_{k-1})^2=\frac{\det{C_k}}{\det{ C_{k-1}}},
\end{equation*}
similarly, for $k+1$:
\begin{equation*}
(a_0)^2*\ldots*(a_{k})^2=\frac{\det{C_{k+1}}}{\det{ C_{k}}}.
\end{equation*}
Two relations above lead to
\begin{equation}
\label{AK} (a_k)^2=\frac{\det{C_{k+1}}\det{
C_{k-1}}}{\det{C_{k}}}, \quad k=1,\ldots, T-1.
\end{equation}
Here we set $\det{C_0}=1,$ $\det{C_{-1}=1}$.

Now using (\ref{matr_eq}) we write down the equation on the last
column of $\left(V^T\right)^{-1}$:
\begin{equation*}
\begin{pmatrix}
{a_{1,1}} & 0 & . &  0 \\
{a_{1,2}} & {a_{2,2}} & 0  &.\\
\cdot & \cdot & \cdot & \cdot  \\
{a_{1,T}} & . & .  & {a_{T,T}}
\end{pmatrix}
\begin{pmatrix}
 c_{1,1} & .. & .. &   c_{1,T} \\
.. & .. & ..  &..\\
\cdot & \cdot & \cdot & \cdot  \\
c_{T,1} &.. &   & c_{T,T}
\end{pmatrix}\begin{pmatrix}
a_{1,T} \\
a_{2,T}\\
\cdot  \\
a_{T,T}
\end{pmatrix}=\begin{pmatrix}
0 \\
0\\
\cdot  \\
1
\end{pmatrix}
\end{equation*}
Note that we know $a_{T,T}$, so we rewrite the above equality in
the form of equation on $(a_{1,T},\ldots,a_{T-1,T})^*$:
\begin{eqnarray*}
\begin{pmatrix}
{a_{1,1}} & 0 & . &  0 \\
{a_{1,2}} & {a_{2,2}} & 0  &.\\
\cdot & \cdot & \cdot & \cdot  \\
{a_{1,T-1}} & . & .  & {a_{T-1,T-1}}
\end{pmatrix}
\begin{pmatrix}
c_{1,1} & .. & .. &   c_{1,T} \\
.. & .. & ..  &..\\
\cdot & \cdot & \cdot & \cdot  \\
 c_{T-1,1} &.. &   &  c_{T-1,T-1}
\end{pmatrix}\begin{pmatrix}
a_{1,T} \\
a_{2,T}\\
\cdot  \\
a_{T-1,T}
\end{pmatrix}\\
+a_{T,T}\begin{pmatrix}
{a_{1,1}} & 0 & . &  0 \\
{a_{1,2}} & {a_{2,2}} & 0  &.\\
\cdot & \cdot & \cdot & \cdot  \\
{a_{1,T-1}} & . & .  & {a_{T-1,T-1}}
\end{pmatrix}\begin{pmatrix}
c_{1,T} \\
c_{2,T}\\
\cdot  \\
c_{T-1,T}
\end{pmatrix}=\begin{pmatrix}
0 \\
0\\
\cdot  \\
0
\end{pmatrix},
\end{eqnarray*}
which is equivalent to the equation
\begin{equation}
\label{matr_eq2}
\begin{pmatrix}
c_{1,1} & .. & .. &  c_{1,T} \\
.. & .. & ..  &..\\
\cdot & \cdot & \cdot & \cdot  \\
 c_{T-1,1} &.. &   &  c_{T-1,T-1}
\end{pmatrix}\begin{pmatrix}
a_{1,T} \\
a_{2,T}\\
\cdot  \\
a_{T-1,T}
\end{pmatrix}=-a_{T,T}\begin{pmatrix}
c_{1,T} \\
c_{2,T}\\
\cdot  \\
c_{T-1,T}
\end{pmatrix}.
\end{equation}
Introduce the notation:
\begin{equation*}
 C_{k-1,k}:=\begin{pmatrix}
 c_{1,1} & .. & .. &   c_{1,k-2} &   c_{1,k}\\
.. & .. & ..  &..\\
\cdot & \cdot & \cdot & \cdot  \\
c_{k-1,1} &.. &   &  c_{k-1,k-2}&  c_{k-1,k}
\end{pmatrix},
\end{equation*}
that is $C_{k-1,k}$ is constructed from $C_{k-1}$ by substituting
the last column by $(c_{1,k},\ldots, c_{k-1,k})^T$. Then from
(\ref{matr_eq2}) we deduce that:
\begin{equation}
\label{a1} a_{T-1,T}=-a_{T,T}\frac{\det{ C_{T-1,T}}}{\det{
C_{T-1}}},
\end{equation}
where we assumed that $\det{C_{-1,0}}=0.$ On the other hand, from
(\ref{AK_form}), (\ref{AK1_form}) we see that
\begin{equation}
\label{a2}
a_{T-1,T}=\left(\prod_{j=0}^{T-1}a_j\right)^{-1}\sum_{k=1}^{T-1}b_k.
\end{equation}
Equating (\ref{a1}) and (\ref{a2}) gives the equalities
\begin{equation*}
\sum_{k=1}^{T-1}b_k=-\frac{\det{C_{T-1,T}}}{\det{C_{T-1}}},\quad
\sum_{k=1}^{T}b_k=-\frac{\det{C_{T,T+1}}}{\det{C_{T}}},
\end{equation*}
from where
\begin{equation}
\label{BK} b_k=-\frac{\det{C_{k,k+1}}}{\det{ C_{k}}}+\frac{\det{
C_{k-1,k}}}{\det{C_{k-1}}}, \quad k=1,\ldots, T-1.
\end{equation}

The relations (\ref{AK}) as in the approach via Krein equations
shows that unlike the case of real Jacobi matrix, the application
of factorization method allows one to recover
$\left(a_k\right)^2$, $k=1,\ldots,T-1$ only. To recover $a_k$ one
need to use the additional information, for example the sequence
of signs. Note that results obtained for dynamic inverse data
corresponds to results obtained to spectral inverse data in
\cite{G}.

We have used two methods to recover the coefficients $(a_k)^2$,
$b_k,$ $k=1,\ldots$. Now we show that impossibility of recovering
$a_k$ is not the weak point of the method, but the feature of the
problem.
\begin{theorem}
\begin{itemize}
\item{1)} For $f\in \mathcal{F}$ the value $u^f_{n,t}$ is odd with
respect to $a_1,a_2,\ldots a_{n-1}$ and even with respect to
$a_n,a_{n+1},\ldots$.
\item{2} The response vector depends on
$\left(a_1\right)^2, \left(a_2\right)^2,\ldots$
\end{itemize}
\end{theorem}
\begin{proof}

The second statement follows from the fact that
$r_{k+1}=u^\delta_{1,k}$, $k=0,1,\ldots$. So we are left to prove
the first one. For $t=1,2,3$ the statement can be checked by
direct calculations, we proceed further by an induction: assuming
that the statement holds for some $t,$ we show it is true for
$t+1$: we write down
\begin{equation}
u_{n,t+1}=-u_{n,t-1}+a_{n}u_{n+1,t}+a_{n-1}u_{n-1,t}+b_nu_{n,t}=0
\end{equation}
The first and the last terms in the right hand side of the above
equality satisfy the statement by an induction assumption. In the
second term the multiple $u_{n+1,t}$ is odd w.r.t. $a_1,\ldots,
a_n$ and even w.r.t. $a_{n+1},\ldots$, as such the second term is
odd w.r.t. $a_1,\ldots, a_{n-1}$ and even w.r.t. $a_{n},\ldots$.
Similarly, the multiple $u_{n-1,t}$ in the third term is odd
w.r.t. $a_1,\ldots, a_{n-2}$ and even w.r.t. $a_{n-1},\ldots$.
Thus after the multiplication by $a_{n-1}$ the third term is odd
w.r.t. $a_1,\ldots, a_{n-1}$ and even w.r.t. $a_n,\ldots$, which
completes the proof.

\end{proof}

\subsection{Characterization of inverse data.}

In the second section we considered the forward problems for
(\ref{Jacobi_dyn}) and (\ref{Jacobi_dyn_aux}): for $(a_0,\dots,
a_{T-1})$, $(b_1,\dots, b_{T-1})$ we constructed the matrices
$W^T$, $W^T_\#$ (\ref{Jac_sol_rep}), (\ref{Goursat}), the response
vector $(r_0,r_1,\dots,r_{2T-2})$ (see (\ref{R_def})) and the
connecting operator  $C_T$ defined in (\ref{C_T_repr}),
(\ref{C_overline_repr}). From Lemma \ref{teor_control} we know
that $C^T$ (and $C_T$) is isomorphism in $\mathcal{F}^T$. We have
also shown that if the vector $(r_0, r_1,\dots,r_{2T-2})$ is a
response vector corresponding to (\ref{Jacobi_dyn}) with the
coefficients $a_0,\dots,a_{T-1}$, $b_1,\dots,b_{T-1},$ then one
can recover those $b_k$ and $\left(a_k\right)^2$ by $a_0=r_0$ and
formulas (\ref{AK}) and (\ref{BK}).

Now we set up a question: can one determine whether a vector
$(r_0,r_1,r_2,\ldots,r_{2T-2})$ is a response vector for dynamical
system (\ref{Jacobi_dyn}) with some $(a_0,\ldots,a_{T-1})$
$(b_1,\ldots,b_{T-1})$? The answer is the following theorem.
\begin{theorem}
\label{Th_char} The  vector $(r_0,r_1,r_2,\ldots,r_{2T-2})$ is a
response vector for the dynamical system (\ref{Jacobi_dyn}) if and
only if the complex symmetric matrix $C^{T-k}$, $k=0,1,\ldots,T-1$
constructed by (\ref{C_T_repr}) is isomorphism in
$\mathcal{F}^{T-k}$.
\end{theorem}
\begin{proof}
First we observe that in the conditions of the theorem we can
substitute $C^T$ by $C_T$ (\ref{C_overline_repr}). The necessary
part of the theorem is proved in the preceding sections. We are
left to prove the sufficiency of this condition.

Let we have a vector $(r_0,r_1,\dots,r_{2T-2})$ such that the
matrix $C_{T}$ constructed from it using (\ref{C_overline_repr})
satisfies conditions of the theorem.

Then we can construct the sequences $(a_0, b_1,\dots,b_{T-1})$
using $a_0=r_0$ and formulas (\ref{BK}) and take arbitrary
sequence of signs to recover $(a_1,\dots,a_{T-1})$ using
(\ref{AK}) and consider the dynamical system (\ref{Jacobi_dyn})
with this coefficients. For this system we construct the response
$(r^{new}_0,r^{new}_1,\dots,r^{new}_{2T-2})$ and connecting
operator $K^{T}$ and its rotated $K_{T}$ using (\ref{C_T_repr})
and (\ref{C_overline_repr}). We will show that the response vector
coincide with the given one.

We have two matrices constructed by (\ref{C_overline_repr}), one
comes from the vector $(r_0,r_1,\dots,r_{2T-2})$ and the other
comes from $(r_0^{new},r^{new}_1,\dots,r^{new}_{2T-2})$. Both of
them have a common property that $C_T$ and $K_T$ are isomorphisms
($C_T$ by the theorem condition and $K_T$ as a connecting
operator). We note that if we calculate the elements of sequences
$(a_1^2,\dots,a_{T-1}^2)$, $(b_1,\dots,b_{T-1})$ using (\ref{AK})
and (\ref{BK}) from any of $C_{T}$ and $K_T$ matrices, we get the
same answer. If so, we obtain that for $k=1,\ldots,T-1$ the
following relations hold:
\begin{eqnarray*}
\frac{\det{C_{k-1,\,k}}}{\det{C_{k-1}}}-\frac{\det{
C_{k,\,k+1}}}{\det{C_{k}}}=\frac{\det{K_{k-1,\,k}}}{{\det{
K}_{k-1}}}-\frac{\det{K_{k,\,k+1}}}{\det{K}_{k}},\\
\frac{\left(\det{C_{k+1}}\right)\left(\det{
C_{k-1}}\right)}{\left(\det{ C_{k}}\right)^2}=\frac{\left(\det{
K_{k+1}}\right)\left(\det{K}_{k-1}\right)}{\left(\det{
K}_{k}\right)^2},\\
\det{C_0}=\det{{K_0}}=1, \quad
\det{C_{-1}}=\det{K_{-1}}=1,\\
\det{C_{-1,\,0}}=\det{K_{-1,\,0}}=0.
\end{eqnarray*}
From these equalities we deduce that
\begin{eqnarray*}
\det{C_k}=\det{K^k},\\
\det{C_{k,\,k+1}}=\det{K_{k,\,k+1}}.
\end{eqnarray*}
The above relations immediately yield that
\begin{equation*}
r_k=r_k^{new}, \quad k=1,\ldots, 2T-2.
\end{equation*}
which completes the proof.

\end{proof}

Below we provide a simple example of the importance of the
condition that the each block $C^{T-k}$, $k=0,1,\ldots T-1$ is an
isomorphism (not only $C^T$ as in the self-adjoint case), cf.
\cite{BH}. We take $r_0=1,$ $r_1=1,$ $r_2=0,$ $r_3=0,$ $r_4=-1$,
in this case in accordance with (\ref{C_T_repr}),
(\ref{C_overline_repr})
\begin{equation}
C_T=\begin{pmatrix} 1 & 1 &0 \\
1 & 1 & 1\\
0 & 1 & 0
\end{pmatrix}
\end{equation}
so $C_T$ is an isomorphism, but $C_{T-1}$ is not invertible and
the formulas (\ref{AK}) and (\ref{BK}) are not applicable.





\begin{thebibliography}{99}

\bibitem{AB}
{S.A. Avdonin, M.I. Belishev,}  \textit{Boundary control and the
dynamic inverse problem for a nonseladjoint Sturm-Louivllle
operator}, {Control and Cybernetics}, Vol. {\bf 25}(3), pp.
429--440, 1996.


\bibitem{AM}
{S.A. Avdonin, V.S. Mikhaylov,} \textit{The boundary control
approach to inverse spectral theory,} {Inverse Problems}, {\bf
26}, no. 4, 045009, 19 pp, 2010.

\bibitem{A}
{F. V. Atkinson} \textit{Discrete and continuous boundary
problems}, Acad. Press, 1964.


\bibitem{B07}
{M.I. Belishev}, {Recent progress in the boundary control method},
\textit{Inverse Problems}, {\bf 23}, R1, 2007.

\bibitem{BH}
{M.I. Belishev, T.Sh. Khabibullin,} {Characterization of data in
dynamical inverse problem for the 1d wave equation with matrix
potential,} Zapiski Nauchnykh Seminarov POMI, , Vol. {\bf 493},
48--72, 2020.

\bibitem{BM_1}
{M.I.Belishev and V.S.Mikhailov}. \textit{Unified approach to
classical equations of inverse problem theory.} {Journal of
Inverse and Ill-Posed Problems}, 20, no 4, 461--488, 2012.

\bibitem{G}
{G. Sh. Guseinov}, \textit{Determination of an infinite
non-self-adjoint Jacobi matrix from its generalized spectral
function,} {Mat. Zametki,} Vol. {\bf 23}(2), 237--248, 1978.



\bibitem{MM1}
{A.S. Mikhaylov, V.S. Mikhaylov}, \textit{Dynamical inverse
problem for the discrete Schr\"odinger operator,} {Nanosystems:
Physics, Chemistry, Mathematics,} Vol. {\bf 7}(5), 842-854, 2016.

\bibitem{MM2}
{A.S. Mikhaylov, V.S. Mikhaylov}, \textit{Dynamic inverse problem
for Jacobi matrices,} {Inverse Problems and Imaging}, Vol {\bf
13}(3), 431-447, 2019.

\bibitem{MM3}
{A.S. Mikhaylov, V.S. Mikhaylov}, \textit{Inverse problem for
dynamical system associated with Jacobi matrices and classical
moment problems,} {Journal of Mathematical Analysis and
Applications}, Vol {\bf 487}(1), 12397, 2020.





\end{thebibliography}
\end{document}